\documentclass[11pt]{article}
\usepackage{epsf}
\usepackage{graphicx}
\usepackage{amsmath}
\usepackage{amssymb}
\usepackage{mathrsfs}
\numberwithin{equation}{section}
\usepackage{srcltx}
\usepackage{color}
\usepackage{todonotes}
\usepackage{hyperref}
\usepackage{background}
\usepackage{upgreek}

\SetBgScale{.1}
\usepackage{color}
\usepackage{amssymb,amsmath,amssymb,amsthm}
\usepackage{hyperref} 
\usepackage{bm} 
\usepackage{mathrsfs}

\newtheorem{theorem}{Theorem}
\newtheorem{definition}{Definition}
\newtheorem{remark}{Remark}

\def\vv{{\bm v}}

\def\yy{{\bm y}}
\def\xx{{\bm x}}
\def\XX{{\bm X}}
\def\nn{{\bm n}}
\def\ff{{\bm f}}

\def\Vdots{\!\mbox{\setlength{\unitlength}{1em}
\begin{picture}(0,0)
\put(-.07,0){.}
\put(-.07,.3){.}
\put(-.07,.6){.}
\end{picture}
}
}

\newcommand\DT[1]{\mathchoice
  {{\buildrel{\hspace*{.1em}\text{\LARGE\bf.}}\over{#1}}}
    {{\buildrel{\hspace*{.1em}\text{\LARGE\bf.}}\over{#1}}}
                 {{\buildrel{\hspace*{.1em}\text{\Large.}}\over{#1}}}
                 {{\buildrel{\hspace*{.1em}\text{\large.}}\over{#1}}}}

\newcommand\pdt[1]{\frac{\partial{#1}}{\partial t}} 
           
\newcommand\R{\mathbb{R}}
\newcommand\N{\mathbb{N}}
\renewcommand\d{{\rm d}}
\newcommand\EE{{\bm e}}
\newcommand{\lineunder}[2]{\LU{\begin{array}[t]{c}\underbrace{#1}\vspace*{.5em}\end{array}}{\mbox{\footnotesize\rm #2}}}
\newcommand{\LU}[2]{\begin{array}[t]{c}#1\vspace*{-1em}\\_{#2}\end{array}}
\newcommand{\linesunder}[3]{\LSU{\begin{array}[t]{c}\underbrace{#1}\vspace*{.5em}\end{array}}{\mbox{\footnotesize\rm #2}}{\mbox{\footnotesize\rm #3}}}
\newcommand{\LSU}[3]{\begin{array}[t]{c}#1\vspace*{-1em}\\_{#2}\vspace*{-.5em}\\_{#3}\end{array}}
\newcommand{\ITEM}[2]{\parbox[t]{.05\textwidth}{{\rm #1}}\hfill\parbox[t]{.95\textwidth}{#2}\vspace*{.8mm}}
\newcommand{\divS}{\mathrm{div}_{\scriptscriptstyle\textrm{\hspace*{-.1em}S}}^{}}

\newcommand{\nablaS}{\nabla_{\scriptscriptstyle\textrm{\hspace*{-.3em}S}}^{}}

\newcommand{\bbI}{\mathbb{I}}

\newcommand{\bbD}{\mathbb D}

\newcommand{\eps}{\varepsilon}

\def\SS{\bm{S}}
\def\FF{\bm{F}}

\begin{document}
\begin{sloppypar}
\allowdisplaybreaks

\centerline{{\LARGE\bf
Visco-elastodynamics at large strains Eulerian}~\footnote{This
research was partially supported also from
the M\v SMT \v CR (Ministry of Education of the Czech Republic) project
CZ.02.1.01/0.0/0.0/15-003/0000493
and the institutional support RVO: 61388998 (\v CR).}}

\bigskip\bigskip

\centerline{\large\sc Tom\'a\v s Roub\'\i\v cek}

\begin{center}
{Mathematical Institute, Charles University,\\
Sokolovsk\'a 83, CZ-186~75~Praha~8, Czech Republic,
\\
and\\
Institute of Thermomechanics of the Czech Academy of Sciences,\\ 
Dolej\v skova 5, CZ-182 00 Praha 8, Czech Republic}

\end{center}

\begin{abstract}
Isothermal visco-elastodynamics in the Kelvin-Voigt rheology is formulated
in the spatial Eulerian coordinates in terms of velocity and deformation
gradient. A generally nonconvex (possibly also frame-indifferent)
stored energy is admitted.
The model involves a nonlinear 2nd-grade nonsimple (multipolar) viscosity
so that the velocity field is well regular. To simplify analytical
arguments,  volume variations of the solid material are
assumed  to be  only rather small so that the mass density is constant,
exploiting the concept of semi-compressible materials.
Existence of weak solutions is proved by using the Galerkin method
combined with a suitable regularization, using
nontrivial results about transport by smooth velocity fields.

\medskip


\medskip


\end{abstract}

\section{Introduction}
Dynamics of deformable elastic or viscoelastic bodies at large strains is one
of basic problems in continuum mechanics. In spite of long-lasting intensive
effort, there still does 
not exist a reasonable analytical theory for existence of reasonably defined
solutions to such problems. Actually, existence of global weak solutions was
articulated as an open problem by J.M.\,Ball in \cite[Problem~12]{Ball02SOPE}
or \cite[Sect.\,1.5]{Ball10PPNE}. There seems to be an agreement that
the basic model for simple purely elastic material at large strains is
not amenable for rigorous analysis and some dissipative mechanisms and
higher gradients (a so-called non-simple material concept) are inevitable.

In the usual notation of {\it deformation} $\yy:\R^d\to\R^d$ considered here
with $d=2$ or 3, one distinguishes the {\it Lagrangian} (referential) {\it
coordinates} $\XX\in\R^d$ and the {\it Eulerian} (i.e.\ actual space)
{\it coordinates} $\xx=\yy(t,\XX)$. In the Lagrangian
coordinates, the deformation formulation of the elastodynamics is
$\varrho\frac{\partial^2}{\partial t^2}\yy={\rm div}_{\!\XX}^{}\SS$
with $\SS=\varphi'(\FF)$ the Piola-Kirchhoff
stress tensor and $\FF=\nabla_{\!\XX}^{}\yy$ is the {\it deformation
gradient}. Restricting on the concept of so-called hyperelastic
materials, the (generalized) Hooke's law $\SS=\varphi'(\FF)$ 
here involves the stored energy $\varphi:\R^{d\times d}\to\R$. This
energy is to be subjected to various requirements, namely frame-indifference
and possibly also a blow-up to $+\infty$ if $\det\FF\to0+$. The latter
requirement will however not be covered by the presented analysis, cf.\
Remarks~\ref{rem-blow-up}.
The former requirement excludes in particular convexity of $\varphi$, although
it complies with polyconvexity, i.e.\ convexity in all subdeterminants of the
deformation gradient. The polyconvexity is well applicable
for static problems but, in the dynamical situations, polyconvex energies
$\varphi$ do not lead to  any  reasonably generalized monotonicity
of  $-{\rm div}\,\varphi'(\nabla\cdot)$ like monotonicity on curl-free
tensor fields, so this concept unfortunately does not seem helpful.

 We will be interested in deformations evolving in time (which are
sometimes called ``motions'').  The formulation in terms of the velocity
and the deformation gradient (briefly $\vv/\FF$ formulation) in the
 Lagrangian  description writes as $\frac{\partial}{\partial t}\FF
=\nabla_{\!\XX}^{}\vv$ and $\varrho\frac{\partial}{\partial t}\vv
={\rm div}_{\!\XX}^{}\SS$ with $\SS=\varphi'(\FF)$. 
It was largely scrutinized in
\cite{Demo00WSCN,DeStTz01VAST,DeStTz12WSUD,KouSpi19QEWS},
exploiting the concept of a so-called measure-valued solution.
In the mere deformation formulation 
$\varrho\frac{\partial^2}{\partial t^2}\yy={\rm div}_{\!\XX}^{}\varphi'(\FF)$,
see also \cite{Rieg03YMSN}. Mostly, the Kelvin-Voigt viscoelasticity
rheology is used and always, various quite restrictive assumptions
must be imposed. In contrast to conventional weak solutions, this concept may,
however, be unacceptably nonselective if not accompanied by some other
attributes, cf.\ the discussion in \cite[Sect.\,8.3]{Roub20ROTV}.
Alternative results concern local-in-time solutions
\cite{DafHru85EMQH,HuKaMa77WPQL}.
Other attempts exploit a stress relaxation, i.e.\
$\SS=\varphi'(\FF)$ is replaced by $\eps\frac{\partial}{\partial t}
(\SS{-}\varphi'(\FF))+(\SS{-}\varphi'(\FF))=0$ with a small
relaxation time $\eps>0$, cf.\ \cite{LatTza06SPSR}, or
a nonlocal (so-called peridynamic) variant of Hooke's law, cf.\
\cite{EmmPhu15MVWS} or the comparison with the conventional local
approach in \cite{EmmPhu15SERN}.

An alternative approach, used mostly for fluids and considered rather
analytically
even more difficult for solids, is to use the Eulerian description, exploiting
the actual deforming configuration, i.e.\ the coordinate
$\xx=\yy(t,\XX)$. Then the velocity reads as $\vv=\DT\xx$ with the dot-notation
$(\cdot)\!\DT{^{}}$ denoting the {\it convective} (also called {\it material})
{\it time derivative}. The chain rule gives the spatial gradient
$\nabla_{\!\xx}^{}\vv=\nabla_{\!\XX}^{}\vv\nabla_{\!\xx}^{}\XX=\DT\FF\FF^{-1}$,
where we used $\FF^{-1}=(\nabla_{\!\XX}^{}\xx)^{-1}=\nabla_{\!\xx}^{}\XX$. As we
will focus here on this Eulerian description, we will omit the subscript $\xx$.
In other words, it gives the {\it transport-and-evolution  equation}
for the deformation-gradient tensor
\begin{align}\label{ultimate}
\DT\FF:=\pdt\FF+(\vv{\cdot}\nabla)\FF=(\nabla\vv)\FF
\end{align}
and the $\vv/\FF$-formulation of the elastodynamics turns into
$\varrho\DT\vv={\rm div}\,{\bm\varSigma}$ with the
Cauchy stress ${\bm\varSigma}$. 
One can also use an Eulerian formulation in terms of $\vv$
and $\FF^{-1}$, as e.g.\ in \cite{CheZha06GESS,SidTho05GETD,Wagn09SHEM}.
The transport-and-evolution rule \eqref{ultimate} was used in incompressible
models with quadratic stored energies in
\cite{CheZha06GESS,HuLin16GSTD,HuWan11GEMD,LeLiZh07GE2D,LiLiZh05HVF,LiuWal01EDFC} or
with convex stored energies in \cite{BFLS18EWSE,KaKoSc21MAWS}, i.e.\ models
at small strains.
 The understanding of \eqref{ultimate} is a bit delicate because it mixes
the Eulerian $\xx$ and the Lagrangian $\XX$; note that $\nabla\vv=\nabla_\xx\vv(\xx)$
while standardly $\FF=\nabla_\XX^{}\yy=\FF(\XX)$. In fact, we consider
$\FF{\circ}\bm\xi$ where $\bm\xi:\xx\mapsto\yy^{-1}(t,\XX)$ is the so-called
{\it return} (sometimes called also a {\it reference}) {\it mapping}. Thus $\FF$
depends on $\xx$ and \eqref{ultimate} is an equality for a.a.\ $\xx$. The
reference mapping $\bm\xi$, which is well defined through its transport equation
$\DT{\bm\xi}=\bm0$, actually does not explicitly occur in the formulation of
the problem. Here we will benefit from the boundary condition $\vv{\cdot}\nn=0$
below, which causes that the actual domain $\varOmega$ does not evolve in time.
The same concerns $\bm{T}$ in \eqref{Euler1+} below, which will make the
problem indeed fully Eulerian, as announced in the title itself.

In general large-strain situations, usually involvement of some viscosity-like
dissipative mechanisms can make analysis more  promising. To this goal,
the {\it Kelvin-Voigt rheology} is most efficient. The elastic stress
(depending on $\FF=\nabla\yy$) is then enhanced by a contribution depending also
on the rate $\nabla\vv$. In the  Lagrangian  description, this
``viscous'' stress is of the form 
$\SS_{\rm v}=\bm\varSigma_{\rm v}(\FF,\nabla\vv)$ with a  function
$\bm\varSigma_{\rm v}:(R^{d\times d})^2\to\R^{d\times d}$ which should
be very nonlinear due to a frame-indifference principle, as pointed out
by S.S.\,Antman \cite{Antm98PUVS}, cf.\ also \cite[Sect.\,9.3]{KruRou19MMCM}
for an analysis. 
In the Eulerian description, incorporation of the Kelvin-Voigt viscosity is
simpler by putting $\SS_{\rm v}=\bm\varSigma_{\rm v}(\EE(\vv))$ with a monotone
$\bm\varSigma_{\rm v}:\R_{\rm sym}^{d\times d}\to\R_{\rm sym}^{d\times d}$  
with $\EE(\vv)=\frac12\nabla\vv^\top\!+\frac12\nabla\vv$.

It is a general understanding that both the  Lagrangian  and the
Eulerian descriptions of both the elastodynamics and visco-elastodynamics at
large strains are not amenable to a reasonable analysis in the sense of mere
global existence of weak solutions, and that involvement of some gradient
theories seems inevitable. In the visco-elastic case, there are naturally
two basic options how to impose some higher-order terms: to the conservative
(elastic) part or to the dissipative (viscous) part, or to both. In the
Lagrangian deformation description, the higher-order deformation gradients 
are analyzed in \cite[Sect.\,9.2-3]{KruRou19MMCM} or \cite{DaRoSt21NHFV}. 
In the Eulerian $\vv/\FF$-formulation, there are even three options
how to involve higher gradients: in addition to  the  two mentioned
options in the Lagrangian variant, one can also regularize the transport
equation \eqref{ultimate} like \eqref{Euler2} or \eqref{Euler2++} below.

The main philosophy of the model below is not to corrupt
the ultimate geometrical relation \eqref{ultimate} governing the transport
of the strain variable $\FF$ by the velocity field $\vv$.
To this goal, we consider a multipolarly dissipative solid
without any stress-diffusion and without strain gradient.
In Section~\ref{sec-model}, we formulate the initial-boundary-value 
for such model and identify formally its energetics. Then, in 
Section~\ref{sec-analysis}, we prove existence of global weak solutions
using the Galerkin method combined with a suitable regularization.
Various modifications and open problems are mentioned, too.

\section{Semi-compressible mass-density-homogeneous\\visco-elastodynamics
}\label{sec-model}

We adopt two assumptions simplifying considerably the analysis and
simultaneously not excluding interesting applications. In particular,
 they weaken  the often imposed incompressibility assumptions
${\rm div}\,\vv=0$ or $\det\FF=1$ with constant mass density $\varrho$, cf.\
Remarks~\ref{rem-mass} and \ref{rem-incompressible}. The incompressible
models, although well applicable in many situations, are physically 
questionable e.g.\ because they do not facilitate propagation 
of pressure (longitudinal) waves. On the other hand, the fully
compressible models are devised rather for 
 gases 
where e.g.\
pressure cannot be negative and zero pressure is related with zero
mass density. In contrast to the incompressible or fully compressible
situations, the {\it semi-compressible} models  (as devised in 
\cite[{Sect.\,5}]{Roub21QISC} for small-strain cases in fluids)
allow for propagation and dispersion also of pressure waves (beside shear
waves), and are suitable for slightly compressible solids  or liquids  without 
 substantial  mass density variations.  This situation
is related with pressures much lower than the elastic bulk modulus. This
modulus is typically quite high in many solid and liquid materials (e.g.\ water
about 2.2\,GPa, magma or rocks about 10\,GPa, steel more than 100\,GPa), but
anyhow considering such materials incompressible (i.e.\ bulk modulus
infinity) would completely suppress effects (like the mentioned pressure-wave
propagation)
which sometimes are of interests.  Simultaneously, the thermomechanical
consistency (as energy balance or frame indifference) is kept  so that
this compromising simplification is an acceptable modelling short-cut
in particular because many analytical technicalities are avoided, cf.\ also
Remark~\ref{rem-mass}. 

We will consider the following visco-elastodynamic system in the
$\vv/\FF$-formulation:
\begin{subequations}\label{Euler+}\begin{align}\nonumber
&\varrho\DT\vv={\rm div}\big({\bm T}+{\bm D}
\big)-\frac\varrho2({\rm div}\,\vv)\,\vv+\ff\ \ \ \text{ with }\ \ {\bm T}=\varphi'(\FF)\FF^\top\!+\varphi(\FF)\bbI
\\[-.3em]
    &\hspace*{16em}
\ \ \text{ and }\ \ {\bm D}=\bbD\EE(\vv)-
{\rm div}\big(\nu|\nabla \EE(\vv)|^{p-2}\nabla\EE(\vv)\big)\,,
\label{Euler1+}
\\\label{Euler2+}
&\DT\FF=(\nabla\vv)\FF
\ \ \ \ \text{ where }\ \ \ \DT\vv=\pdt\vv+(\vv{\cdot}\nabla)\vv\ \ \text{ and }\ \ 
\DT\FF=\pdt\FF+(\vv{\cdot}\nabla)\FF\,,
 \end{align}\end{subequations}
 where ${\bm T}$ is the conservative part of the
Cauchy stress, while ${\bm D}$ is its dissipative part which
contains the viscous-moduli tensor $\bbD$ and
also the contribution of the so-called hyperstress
$\nu|\nabla \EE(\vv)|^{p-2}\nabla\EE(\vv)$ with some
coefficient $\nu>0$ assumed to be small.  The term
$-\frac\varrho2({\rm div}\,\vv)\,\vv$ in the first equation \eqref{Euler1+} is
a force introduced by Temam \cite{Tema69ASEN} rather for numerical
purposes to balance the energy. Actually, it violates the Galilean invariance,
as pointed out in \cite{Toma21ITST} where a certain justification can be found.
Yet, this Temam's force is presumabl{y}  very small in media which
are only very slightly compressible, so-called quasi-incompressible.
This is the price for a simplification that the mass density $\varrho$ is
constant, cf.\ Remark~\ref{rem-mass} for a full model.
In \eqref{Euler+}, the notation ``$\,\cdot\,$'' denotes the scalar products
of vectors and later ``$\,:\,$'' and
``$\,\Vdots\,$'' will be used for the scalar products of 
matrices and 3rd-order tensors, respectively.

Together with the standard linear viscous stress
$\bbD\EE(\vv)$  in \ref{Euler1+}, we have employed the concept of 
the so-called {\it nonsimple fluids},
  devised by E.\,Fried and M.\,Gurtin \cite{FriGur06TBBC}
and earlier, even more generally and nonlinearly as multipolar fluids, by
J.\,Ne\v cas at al.\ \cite{BeBlNe92PBMV,NeNoSi89GSIC,NecRuz92GSIV}
or solids \cite{Ruzi92MPTM,Silh92MVMS}. More specifically, we use 2nd-grade
nonsimple fluids, also called {\it bipolar fluids}, in a nonlinear variant.
 Here it leads to the viscous {\it hyperstress}
$\nu|\nabla \EE(\vv)|^{p-2}\nabla\EE(\vv)$; the
preposition ``hyper'' means that  it contributes to the stress through its
 divergence. Here, the goal is to ensure
$\nabla\vv\in L_{\rm w*}^1(I;L^\infty(\varOmega;\R^{d\times d}))$, cf.\ 
in particular the estimates \eqref{test-FF}, \eqref{test-vv}, or \eqref{test-Delta-FF} below.

We will consider a fixed bounded Lipschitz domain $\varOmega$ in $\R^d$
with the boundary $\varGamma$. Then, the system  of two nonlinear parabolic
equations \eqref{Euler+} should be accompanied by some boundary conditions,
e.g.\
\begin{align}\label{BC}
  &
\vv{\cdot}\nn=0\,,\ \ \ \ \big[  ({\bm T}{+}{\bm D})\nn{+}\divS(
\nu|\nabla\EE(\vv)|^{p-2}(\nabla\EE(\vv))\nn)\big]_\text{\sc t}^{}=\bm g\,,
   \ \ \ \ \nabla\EE(\vv){:}(\nn{\otimes}\nn)=0\,,
\end{align}
with $\nn$ denoting the unit outward normal to $\varGamma$,
$[\cdot]_\text{\sc t}^{}$ the tangential component of a vector, i.e.\
$[\bm\sigma]_\text{\sc t}^{}=\bm\sigma-(\bm\sigma{\cdot}\nn)\nn$
for a vector $\bm\sigma$. 
In \eqref{BC}, $\divS={\rm tr}(\nablaS)$ denotes the $(d{-}1)$-dimensional
surface divergence with ${\rm tr}(\cdot)$ being the trace of a
$(d{-}1){\times}(d{-}1)$-matrix and
$\nablaS v=\nabla v-\frac{\partial v}{\partial\nn}\nn$ 
being the surface gradient of $v$.
Naturally, ${\bm g}{\cdot}\nn=0$ is to be assumed
 if we want to recover the boundary conditions \eqref{BC} in the
classical form, otherwise the weak form does not directly need it.

To reveal at least formally the energetics behind the system \eqref{Euler+},
one should test \eqref{Euler1+} by $\vv$ and \eqref{Euler2+} by $\SS$.
We use several calculations exploiting integration over $\varOmega$ and
Green's formula.

For the convective term in \eqref{Euler1+} we use the calculus
\begin{align}\nonumber
  &\int_\varOmega\varrho(\vv{\cdot}\nabla)\vv{\cdot}\vv\,\d x
 =\int_\varGamma\varrho|\vv|^2(\vv{\cdot}\nn)\,\d S
  - \int_\varOmega\vv\nabla(\varrho\vv\otimes\vv)\,\d x
  \\&\qquad\nonumber=\int_\varGamma\varrho|\vv|^2(\vv{\cdot}\nn)\,\d S
  - \int_\varOmega\varrho|\vv|^2{\rm div}\,\vv
  +\varrho\vv{\cdot}\nabla\vv\cdot\vv
+|\vv|^2(\nabla\varrho{\cdot}\vv)\,\d x
  \\&\qquad=\int_\varGamma\frac\varrho2|\vv|^2(\vv{\cdot}\nn)\,\d S
  -\int_\varOmega\frac\varrho2|\vv|^2({\rm div}\,\vv)
  +\frac12|\vv|^2(\nabla\varrho{\cdot}\vv)\,\d x
  =-\!\int_\varOmega\frac\varrho2|\vv|^2({\rm div}\,\vv)\,\d x
  \label{convective-tested}\end{align}
which also uses the first boundary condition \eqref{BC} and the assumption that
$\varrho$  is  constant so that $\nabla\varrho=0$.
Here the role of Temam's bulk force
$\frac12\varrho({\rm div}\,\vv)\vv$ in \eqref{Euler1+} is revealed.

Furthermore, we use the algebra  ${\bm S}{:}({\bm V}\FF)=({\bm S}\FF^\top){:}{\bm V}$ for ${\bm V}=\nabla\vv$ and ${\bm S}=\varphi'(\FF)$ 
and the following calculus
\begin{align}\nonumber
\int_\varOmega{\bm T}{:}\nabla\vv\,\d x
&=\int_\varOmega\big(\varphi'(\FF)\FF^{\top}\!\!+\varphi(\FF)\bbI\big){:}\nabla\vv\,\d x=\int_\varOmega(\nabla\vv)\FF^{\top}\!{:}\varphi'(\FF)
+\varphi(\FF){\rm div}\,\vv\,\d x
\\&\nonumber=\int_\varOmega\DT\FF{:}\varphi'(\FF)
+\varphi(\FF){\rm div}\,\vv\,\d x
=\int_\varOmega\varphi'(\FF){:}\pdt{\FF}+\hspace{-.9em}\lineunder{\varphi'(\FF){:}(\vv{\cdot}\nabla)\FF}{$=\nabla\varphi(\FF){\cdot}\vv$}\hspace{-.9em}
+\varphi(\FF){\rm div}\,\vv\,\d x
\\[-.9em]&=\frac{\d}{\d t}\int_\varOmega\varphi(\FF)\,\d x+\int_\varGamma\varphi(\FF)\hspace{-.8em}\lineunder{\vv{\cdot}\nn}{$=0$}\hspace{-.8em}\,\d S\,,
\label{Sstr-calculus}\end{align}
where we used also the Green formula  for
$$
\int_\varOmega\nabla\varphi(\FF){\cdot}\vv+\varphi(\FF){\rm div}\,\vv\,\d x
=\int_\varGamma{\rm div}(\varphi(\FF)\vv)\,\d x
=\int_\varGamma\varphi(\FF)\vv{\cdot}\nn\,\d S=0\,,
$$
exploiting the  boundary condition $\vv{\cdot}\nn=0$.
For the term ${\rm div}^2(|\nabla\EE(\vv)|^{p-2}\nabla\EE(\vv))$ tested
by $\vv$, we use twice Green{'s}  formula together with a surface Green
formula and the boundary conditions \eqref{BC}, cf.\ also
\cite[Sect.\,2.4.4]{Roub13NPDE} for technical details. This procedure yields
(at least formally) the energy balance:
\begin{align}
  &\frac{\d}{\d t}
  \int_\varOmega\!\!\!\!\linesunder{\frac\varrho2|\vv|^2}{kinetic}{energy}\!\!\!\!\!\!+\!\!\!\!\!\!\linesunder{\varphi(\FF)\!_{_{_{_{_{_{}}}}}}}{stored}{energy}\!\!\!\!\!\,\d x
+\int_\varOmega\!\!\!\!\lineunder{
\bbD \EE(\vv){:}\EE(\vv)+\nu|\nabla \EE(\vv)|^p_{_{_{}}}
\!}{dissipation rate}\!\!\!\d x=\int_\varOmega\!\!\!\!\!\!\!\!\!\!\linesunder{\ff{\cdot}\vv_{_{_{_{_{}}}}}\!\!\!}{power of}{external load}\!\!\!\!\!\!\!\!\!\d x+\int_\varGamma\!\!\!\!\!\!\!\!\!\!\linesunder{\bm g{\cdot}\vv_{_{_{_{_{}}}}}}{power of}{traction load}\!\!\!\!\!\!\!\!\!\!\!\d S\,.
  \label{energy+++}\end{align}

Further, we prescribe the initial conditions, i.e.
\begin{align}\label{IC}
\vv|_{t=0}^{}=\vv_0\ \ \ \text{ and }\ \ \ \FF|_{t=0}^{}=\FF_0\,.
\end{align}

\begin{remark}[{\sl Varying mass density}]\label{rem-mass}\upshape
   The above model used the simplification based on the assumption of a
 constant mass density $\varrho$, cf.\ the calculus \eqref{convective-tested},
   and neglects its variations during
   volumetric deformation (which is typically indeed small in
   liquids and in solids, too, in contrast to gases).
In fact, the initial (possibly spatially nonconstant) density should evolve in
time by the continuity equation $\DT\varrho=-\varrho\,{\rm div}\,\vv$. In the
context of solids in Eulerian description, see also
\cite{GuFrAn10MTC,HuWan11GEMD,LiuWal01EDFC,Mart19PCM,SidTho05GETD}.
Having omitted this continuity equation has
been here compensated by the 
 extra  force $-\varrho({\rm div}\,\vv)\vv/2$
in \eqref{Euler1+}, which is presumably very small.
This simplifies a lot of calculations and analytical arguments. Here,
the analysis would likely be doable
for the full model involving also the continuity equation  with an
initial condition $\varrho|_{t=0}^{}=\varrho_0$ or,
equivalently, the algebraic relation $\varrho=\varrho_0/\det\FF$ 
and omitting the extra compensating  force
because the concept of nonlinear nonsimple material
ensures $\nabla\vv$
in $L_{\rm w*}^p(I;L^\infty(\varOmega;\R^{d\times d}))$ which in turn allows
for estimation of $\nabla\varrho$, cf.\
 also Remark~\ref{rem-blow-up} below.
\end{remark}

\begin{remark}[{\sl The Cauchy stress ${\bm T}$ alternatively}]\upshape
The pressure contribution to the Cauchy stress ${\bm T}$ in \eqref{Euler1+}
is related to that the stored energy $\varphi$ is counted truly Eulerian, i.e.\
per actual volume. Another approach is to count the stored energy per the
referential volume (let us denote it by $\varphi_\text{\sc r}$).
The relation between $\varphi$ and $\varphi_\text{\sc r}$ is
$\varphi(\FF)=\varphi_\text{\sc r}(\FF)/\det\FF$.
In terms of $\varphi_\text{\sc r}$, the conservative part of the Cauchy stress
${\bm T}$ transforms to
\begin{align}\nonumber
{\bm T}=\varphi'(\FF)\FF^\top\!\!+\varphi(\FF)\bbI
&=\bigg(\frac{\varphi_\text{\sc r}'(\FF)}{\det\FF}
-\frac{\varphi_\text{\sc r}(\FF){\rm Cof}\FF}{(\det\FF)^2}\bigg)\FF^\top+
\frac{\varphi_\text{\sc r}(\FF)}{\det\FF}\bbI
\\&
=\ \frac{\varphi_\text{\sc r}'(\FF)-\varphi_\text{\sc r}(\FF)\FF^{-\top}\!\!\!\!}
{\det\FF}\FF^\top\!+
\frac{\varphi_\text{\sc r}(\FF)}{\det\FF}\bbI
\ =\ \frac{\varphi_\text{\sc r}'(\FF)}{\det\FF}\FF^\top\,.
\end{align}
 Then 
$\varphi(\FF)$ in \eqref{energy+++} should be replaced by
$\varphi_\text{\sc r}(\FF)/\det\FF$, cf.\ e.g.\ \cite{GiKiLi17VMCF}
or \cite[Ch.7]{Gurt83TFE}. 
This variant would need to control $\det\FF$, which we avoided in this paper,
cf.\ also Remark~\ref{rem-blow-up} below for a possible modification.
\end{remark}

\begin{remark}[{\sl Conservative multipolar variant with stress
diffusion}]\label{rem-diffusion}\upshape
The multipolar concept can alternatively use a conservative higher gradient
instead of the dissipative higher-order term ${\rm div}
(\nu|\nabla\EE(\vv)|^{p-2}(\nabla\EE(\vv))$
in \eqref{Euler1+}, as often occurs in literature since the work by
R.A.\,Toupin \cite{Toup62EMCS} and R.D.~Mindlin \cite{Mind64MSLE}, cf.\ e.g.\
\cite{FosRoy02LMHC,PoGiVi10HHCS}. Yet, this is more complicated because such
term is reflected in  a nonlinear  (Korteweg-like) symmetric 
 capillarity  stress $\bm{K}$. In addition, for analytical reasons, the
 transport-and-evolution equation
\eqref{Euler2+} for $\FF$ must be enhanced by stress diffusion. 
This leads to the visco-elastodynamic system
 \begin{subequations}\label{Euler}\begin{align}\nonumber
&\varrho\Big(\pdt\vv+(\vv{\cdot}\nabla)\vv\Big)
={\rm div}\big(\bm{T}+\bm{D}
+\bm{K}\big)-\frac\varrho2({\rm div}\,\vv)\,\vv+\ff\,,
\ \ \ \ \text{where }\ 
\bm{T}=\SS\FF^\top\!\!+\varphi(\FF)\bbI\,,\ \\[-.3em]
&\hspace*{6.5em}
\ \ 
\SS=\varphi'(\FF)-\varkappa\Delta\FF\,,
\ \ \ \ 
    \label{Euler1}
\bm{D}=\bbD\EE(\vv)\,,
 \ \ \text{ and }\ \ \:
 \bm{K}=\frac\varkappa2|\nabla\FF|^2\bbI-\varkappa
 \nabla\FF{\otimes}\nabla\FF\,,
\\[-.3em]\label{Euler2}
&\pdt\FF+(\vv{\cdot}\nabla)\FF=(\nabla\vv)\FF+\epsilon\Delta\SS\
\end{align}\end{subequations}
with the boundary conditions $\vv{\cdot}\nn=0$,
$(\vv{\cdot}\nabla)\FF=0$, and $(\vv{\cdot}\nabla)\SS=0$.
 The capillarity coefficient $\varkappa>0$ and a stress-diffusion
coefficient $\epsilon>0$ are presumably small. 
Corrupting the ``geometrical''  transport-and-evolution 
rule \eqref{ultimate} in 
\eqref{Euler2} might be rather controversial and purely mathematically
motivated, as actually openly articulated in
\cite{BFLS18EWSE,KaKoSc21MAWS} for the variant of $\Delta\FF$
instead of $\Delta\SS$. Anyhow, in fluid dynamics,
 the stress diffusion $\epsilon\Delta\SS$
 was advocated in series of works by H.\,Brenner,
 cf.\ e.g.\ \cite{Bren05KVT,Bren06FMR}.
 Cf.\ also a discussion in \cite{OtStLi09IDCM} and a thermodynamical
 justification in \cite{VaPaGr17EMFF}. Here, \eqref{Euler2} with
$\SS=\varphi'(\FF)-\varkappa\Delta\FF$ is formally rather Cahn-Hilliard's
diffusion of $\FF$.
The energy balance is again obtained by the test of \eqref{Euler1} by $\vv$
and of \eqref{Euler2} by $\SS$. Using the calculus
\eqref{test-Delta} below, we obtain the energy balance
\begin{align}
  &\frac{\d}{\d t}
  \int_\varOmega\!\!\!\!\linesunder{\frac\varrho2|\vv|^2}{kinetic}{energy}\!\!\!\!\!\!+\!\!\!\!\!\linesunder{\varphi(\FF)+\frac{\varkappa}2|\nabla\FF|^2}{stored}{energy}\!\!\!\!\d x
  +\int_\varOmega\!\!\!\!\lineunder{\bbD \EE(\vv){:}\EE(\vv)
+\epsilon|\nabla\SS|^2_{_{_{}}}\!}{dissipation rate}\!\!\!\d x=\int_\varOmega\!\!\!\!\!\!\!\!\!\!\!\linesunder{\ff{\cdot}\vv_{_{_{_{_{}}}}}\!\!\!}{power of}{external load}\!\!\!\!\!\!\!\!\!\!\d x+\int_\varGamma\!\!\!\!\!\!\!\!\!\!\linesunder{\bm g{\cdot}\vv_{_{_{_{_{}}}}}}{power of}{traction load}\!\!\!\!\!\!\!\!\!\!\!\d S\,.
  \label{energy++}\end{align}
From this, one can read the estimate 
$\SS\in L^2(I;H^1(\varOmega;\R^{d\times d\times d}))$. Then, assuming
$\varphi''$ bounded, one has $\varkappa\Delta(\nabla\FF)
=\nabla(\varphi'(\FF){-}\SS)
=\varphi''(\FF)\nabla\FF-\nabla\SS\in L^2(I{\times}\varOmega;\R^{d\times d\times d})$
and, assuming $\varOmega$ smooth, one can further use the
$H^2$-regularity for $\varkappa\Delta$, so that
$\nabla\FF\in L^2(I;H^2(\varOmega;\R^{d\times d\times d}))$.
By comparison from \eqref{Euler2}, 
$\pdt{}\FF=(\nabla\vv)\FF-(\vv{\cdot}\nabla)\FF+\epsilon\Delta\SS$
belongs to $L^{2}(I;L^1(\varOmega;\R^{d\times d\times d})
{+}H^1(\varOmega;\R^{d\times d\times d})^*)$ so that,
by Aubin-Lions' theorem,
$\nabla\FF$ and $\nabla\FF{\otimes}\nabla\FF$ are ``compact'',
which is needed for the capillarity-type  stress  $\bm{K}$.
Such alternative gradient model  in some sense  improves some previous
strain-diffusion models \cite{BFLS18EWSE,KaKoSc21MAWS}.
\end{remark}

\section{Analysis of the model (\ref{Euler+})}\label{sec-analysis}

We now prove existence of weak solutions by a constructive
method, i.e.\ by a suitable approximation and its convergence.
This is rather technical in general. The time discretisation (Rothe's method)
standardly needs convexity of $\varphi$ (which is not a realistic assumption
here) possibly weakened if there is some viscosity in $\FF$ (which is not
directly considered here, however). The  conformal  space
discretisation ({the Faedo-Galerkin}  method) cannot
directly copy the energetics because the test
of \eqref{Euler2+} by $\SS$ is problematic in  this 
approximation
as $\SS=\varphi'(\FF)$ is not in the respective finite-dimensional
space in general. We facilitate the analytical issue here by imposing
rather strong growth assumption \eqref{ass1} on $\varphi$ below, which might be 
quite restrictive, cf.\ Remark~\ref{rem-blow-up}.

We will use the standard notation concerning the Lebesgue and the Sobolev
spaces  on the domain $\varOmega\subset\R^d$, as actually already employed in Remark~\ref{rem-diffusion}. Namely,  for $n\in\N$, 
$L^p(\varOmega;\R^n)$ denotes the Banach spaces of Lebesgue measurable functions
$\varOmega\to\R^n$ whose Euclidean norm is integrable with $p$-power, and
$W^{k,p}(\varOmega;\R^n)$ the space of functions from $L^p(\varOmega;\R^n)$ whose
derivatives of the order $k$ are in $L^p(\varOmega;\R^{n\times kd})$.
 If $n=1$, we will write simply $L^p(\varOmega)$ or $W^{k,p}(\varOmega)$. 
Moreover, $W_0^{2,p}(\varOmega;\R^d):=\{\vv\in W^{2,p}(\varOmega;\R^d);\
\vv{\cdot}\nn=0\text{ on }\varGamma\}$.
We also write briefly $H^k=W^{k,2}$. Moreover, for a Banach space
$X$ and for $I=[0,T]$, we will use the notation $L^p(I;X)$ for the Bochner
space of Bochner measurable functions $I\to X$ whose norm is in $L^p(I)$, 
and $H^1(I;X)$ for functions $I\to X$ whose distributional derivative
is in $L^2(I;X)$. Moreover, $(\cdot)^*$ will denote the dual space and
$p'=p/(p{-}1)$ is the conjugate exponent with the convention
$p'=\infty$ for $p=1$.
Occasionally, we will use $L_{\rm w*}^p(I;X)$ for weakly* measurable
functions $I\to X$ for nonseparable spaces $X$ which are duals to some other
Banach spaces (specifically for $L^\infty(\varOmega)$).

Let us first summarize the assumptions:
\begin{subequations}\label{ass}
\begin{align}\nonumber
&\varphi:\R^{d\times d}
\to[0,+\infty]\ \text{ continuously differentiable and }\ 
\\
&\qquad\qquad\exists\,
\ell\in\R
\ \forall F\in\R^{d\times d}:\ \ 
0\le\varphi(F)\le \ell(1+|F|)\ \ \ \text{ and }\ \ \ |\varphi'(F)|\le \ell\,,
\label{ass1}\\&
\ff\in L^1(I;L^2(\varOmega;\R^d))+L^2(I;L^1(\varOmega;\R^d)),\ \ \
{\bm g}\in L^2( I;L^1(\varGamma;\R^d)),
\label{ass2}\\&\vv_0\in L^2(\varOmega;\R^d),\ \ \ \FF_0\in 
H^1(\varOmega;\R^{d\times d}),\ \ \
\label{ass3}\\[-.3em]&\varrho,\nu>0,\ \ \ \bbD\in(\R_{\rm sym}^{d\times d})^2
\ \text{ positive definite},\ \ p>d\,.
\label{ass4}\end{align}\end{subequations}
Beside, it is also reasonable
to assume the frame indifference, i.e.\ $\varphi(F)=\varphi(QF)$
for any $Q\in{\rm SO}(d)=\{A\in\R^{d\times d};\ A^\top\!A=AA^\top\!=\bbI,\ \det A=1\}$.  It is equivalent to $\varphi(F)=\phi(F^\top F)$ for some
$\phi$ so that ${\bm T}=\varphi'(\FF)\FF^\top=2\FF\phi'(\FF^\top\FF)\FF^\top$
and thus the whole Cauchy stress ${\bm T}+{\bm D}$ is symmetric. 
Yet, we will not explicitly need it for our analysis. 

For the definition of a weak solution, we do not need to impose
a-priori smoothness of $\FF$ when using the calculus
$\int_\varOmega(\vv{\cdot}\nabla)\FF{:}\widetilde\FF\,\d x=
\int_\varGamma\FF{:}\widetilde\FF(\vv{\cdot}\nn)\,\d S-
\int_\varOmega\FF{:}(({\rm div}\,\vv)\widetilde\FF+
(\vv{\cdot}\nabla)\widetilde\FF)\,\d x$, although later we will prove even
quite high regularity of $\FF$. Like already in Section~\ref{sec-model},
for the term ${\rm div}^2(|\nabla\EE(\vv)|^{p-2}\nabla\EE(\vv))$ tested
by $\widetilde\vv$ with $\nn{\cdot}\widetilde\vv=0$, we will use twice the Green
formula on $\varOmega$ together with a surface Green
formula on $\varGamma$ together with the boundary conditions \eqref{BC}.
In addition, we use the by-part integration in time.
This results to:

\medskip

\begin{definition}[Weak solutions to \eqref{Euler+}--\eqref{BC}]\label{def}
A couple $(\vv,\FF)$ with $\vv\in L^p(I;W^{2,p}(\varOmega;\R^d))\cap
H^1(I;L^2(\varOmega;\R^d))$
with $\vv{\cdot}\nn=0$ on $I{\times}\varGamma$
and $\FF\in L^2(I{\times}\varOmega;\R^{d\times d})$ 
is called a weak solution to the initial-boundary-value
problem \eqref{Euler+}--\eqref{BC} with \eqref{IC} if the
integral identity 
\begin{subequations}\begin{align}\nonumber
&\int_0^T\!\!\!\int_\varOmega\bigg(\varphi'(\FF)\FF^\top{:}
\nabla\widetilde\vv
+\Big(\frac\varrho2({\rm div}\,\vv)\vv
+\varrho(\vv{\cdot}\nabla)\vv\Big){\cdot}\widetilde\vv
-\varrho\vv{\cdot}\frac{\partial\widetilde\vv}{\partial t}
-\varphi(\FF){\rm div}\,\widetilde\vv
\\[-.4em]&\quad
+\nu|\nabla\EE(\vv)|^{p-2}\nabla\EE(\vv)\Vdots\nabla\EE(\widetilde\vv)
\bigg)\,\d x\d t
=\int_0^T\!\!\!\int_\varOmega\!\ff{\cdot}\widetilde\vv\,\d x\d t
+\int_0^T\!\!\!\int_{\varGamma}\bm g{\cdot}\widetilde\vv\,\d x\d t
+\int_\varOmega\varrho\vv_0{\cdot}\widetilde\vv(0)\,\d x
\label{def1}
\intertext{holds for any $\widetilde\vv\in L^p(I;W^{2,p}(\varOmega;\R^d))\cap
H^1(I;L^2(\varOmega;\R^d))$ with $\widetilde\vv(T)=0$ and
$\widetilde\vv{\cdot}\nn=0$ on $I{\times}\varGamma$, and also the
integral identity}
&\int_0^T\!\!\!\int_\varOmega
(\nabla\vv)\FF{:}\widetilde\FF+({\rm div}\,\vv)\FF{:}\widetilde\FF
+\FF{:}(\vv{\cdot}\nabla)\widetilde\FF
+\FF{:}\frac{\partial\widetilde\FF}{\partial t}\,\d x\d t
+\int_\varOmega\!\FF_0{:}\widetilde\FF(0)\,\d x=0
\label{def2}
\end{align}\end{subequations}
holds for any $\widetilde\FF\in W^{1,\infty}(I{\times}\varOmega;\R^{d\times d})$ with
$\widetilde\FF(T)=0$.
\end{definition}

\begin{theorem}[Existence of weak solutions]
Let \eqref{ass} hold. Then:\\
\ITEM{(i)}{There exists a weak solution $(\vv,\FF)$ according  to 
Definition~\ref{def} such that also
$\FF\in L^\infty(I;H^1(\varOmega;\R^{d\times d}))$, the transport-and-evolution
rule \eqref{Euler2+} holds a.e.\ on $I{\times}\varOmega$, and also the energy
balance \eqref{energy+++}  integrated over a time interval $[0,t]$ is satisfied
for all $t\in I$.}
\ITEM{(ii)}{If $\varphi$ is also twice continuously
differentiable with $\varphi''$ bounded, then
$\SS=\varphi'(\FF)\in L^\infty(I;H^1(\varOmega;\R^{d\times d\times d}))$.}
\end{theorem}

\begin{proof}
For  clarity, we divide the proof into four steps.

\medskip{\it Step 1: Galerkin approximation of a regularized problem}.
Instead of \eqref{Euler2+}, we will consider in this step a regularization 
as in \cite{BFLS18EWSE,KaKoSc21MAWS}:
\begin{align}
&\pdt\FF+(\vv{\cdot}\nabla)\FF=(\nabla\vv)\FF+\eps\Delta\FF
\label{Euler2++}
\end{align} 
with the boundary conditions \eqref{BC} enhanced also by other $d{\times}d$
conditions $(\nn{\cdot}\nabla)\FF=0$.
The weak solution is defined analogously as in Definition~\ref{def} but
with \eqref{def2} enhanced by the term
$-\eps\nabla\FF\Vdots\nabla\widetilde\FF$.

We will then use the Faedo-Galerkin method exploiting a
nested sequence of finite-dimensional sequences
$V_k\subset W_{0}^{2,\infty}(\varOmega)$
such that $\bigcup_{k\in\N}V_k$ is dense in
$W_0^{2,p}(\varOmega)$ and composed with eigenfunctions of the self-adjoint
$\Delta$-operator with the homogeneous Neumann boundary condition;
recall the notation $W_0^{2,p}(\varOmega)$ used here for vectorial 
functions with normal boundary traces zero.
This ensures that $\Delta\vv\in V_k$ provided $\vv\in V_k$.
We also approximate the initial conditions: $\vv_{0,k}\in V_k^d$
and $\FF_{0,k}\in V_k^{d\times d}$ so that $\{\vv_{0,k}\}_{k\in\N}$ is bounded in
$L^2(\varOmega;\R^d)$ and $\{\FF_{0,k}\}_{k\in\N}$ is bounded in
$H^1(\varOmega;\R^{d\times d})$ and $\vv_{0,k}\to\vv_0$ and
$\FF_{0,k}\to\FF_0$ strongly in these spaces.

Then we consider the Galerkin approximate solution
$(\vv_{\eps k},\FF_{\eps k}):I\to V_k^d\times V_k^{d\times d}$ to the regularized
system \eqref{Euler1+} and \eqref{Euler2++}. Existence of such solution is
based by standard ordinary-differential-equation arguments with a successive
prolongation exploiting the $L^\infty(I)$-estimates below.

\medskip{\it Step 2: estimates for the approximated regularized problem}.
To derive basic a-priori estimates for such Galerkin approximation,
we test \eqref{Euler2++} in its Galerkin approximation by $\FF_{\eps k}$. Thus
we obtain the estimate
\begin{align}\nonumber
&\eps\int_\varOmega|\nabla\FF_{\eps k}|^2\,\d x
+\frac{\d}{\d t}\int_\varOmega\frac12|\FF_{\eps k}|^2\,\d x
=\int_\varOmega\big((\nabla\vv_{\eps k})\FF_{\eps k}-(\vv_{\eps k}{\cdot}\nabla)\FF_{\eps k}\big){:}\FF_{\eps k}\,\d x
\\[-.4em]&\quad
=\int_\varOmega(\nabla\vv_{\eps k})\FF_{\eps k}{:}\FF_{\eps k}+\frac12({\rm div}\,\vv_{\eps k})|\FF_{\eps k}|^2\,\d x
\le
\frac32\|\nabla\vv_{\eps k}\|_{L^\infty(\varOmega;\R^{d\times d})}^{}\|
\FF_{\eps k}\|_{L^2(\varOmega;\R^{d\times d})}^2\,;
\label{test-FF}\end{align}
here we used also the calculus
\begin{align}
  \int_\varOmega(\vv{\cdot}\nabla)\FF{:}\FF\,\d x
  &=\int_\varGamma|\FF|^2(\vv{\cdot}\nn)\,\d S
  -\int_\varOmega\!\FF{:}(\vv{\cdot}\nabla)\FF+({\rm div}\,\vv)|\FF|^2\,\d x
=-\frac12\int_\varOmega({\rm div}\,\vv)|\FF|^2\,\d x
\label{calcul-FF}\end{align}
together with the boundary condition $\vv{\cdot}\nn=0$. Moreover,
testing \eqref{Euler1+}  in its Galerkin approximation by $\vv_{\eps k}$,
we obtain
\begin{align}\nonumber
&\int_\varOmega\bbD\EE(\vv_{\eps k}){:}\EE(\vv_{\eps k})+\nu|\nabla\EE(\vv_{\eps k})|^p\,\d x
+\frac{\d}{\d t}\int_\varOmega\frac\varrho2|\vv_{\eps k}|^2\,\d x
\\[-.4em]&\nonumber\qquad
\le\int_\varOmega\!\varphi'(\FF_{\eps k})\FF_{\eps k}^\top{:}\nabla\vv_{\eps k}
+\varphi(\FF_{\eps k}){\rm div}\,\vv_{\eps k}+\ff{\cdot}\vv_{\eps k}\,\d x
+\!\int_\varGamma\bm g{\cdot}\vv_{\eps k}\,\d S
\\[-.0em]&\nonumber\qquad\le
\int_\varOmega 
2\ell(1+|\FF_{\eps k}|)|\nabla\vv_{\eps k}|+\ff{\cdot}\vv_{\eps k}\,\d x
+\!\int_\varGamma\bm g{\cdot}\vv_{\eps k}\,\d S
\\[-.1em]&\nonumber\qquad\le
C_{\ell,\delta}\|\FF_{\eps k}\|_{L^2(\varOmega;\R^{d\times d})}^2
+\delta\|\nabla\vv_{\eps k}\|_{L^2(\varOmega;\R^{d\times d})}^2
+\|\ff_1\|_{L^2(\varOmega;\R^d)}^2\big(1+\|\vv_{\eps k}\|_{L^2(\varOmega;\R^d)}^2\big)
\\[-.1em]&\hspace*{5em}
+\big(\|\ff_2\|_{L^1(\varOmega;\R^d)}^2\!+\|\bm g\|_{L^1(\varGamma;\R^d)}^2\big)/\delta
+\delta\|\vv_{\eps k}\|_{L^\infty(\varOmega;\R^d)}^2
+\delta\|\vv_{\eps k}|_\varGamma\|_{L^\infty(\varGamma;\R^d)}^2\,,
\label{test-vv}\end{align}
where $\ell$ is from \eqref{ass1} and where $\ff=\ff_1+\ff_2$ with some
$\ff_1\in L^1(I;L^2(\varOmega;\R^d))$ and $\ff_2\in L^2(I;L^1(\varOmega;\R^d))$
referring to \eqref{ass2}. The last terms in \eqref{test-vv} are
still to be estimated using Korn's inequality
and the boundedness of the embedding/trace operators for
$\max(\|\vv\|_{L^\infty(\varOmega;\R^d)},\|\vv|_\varGamma\|_{L^\infty(\varGamma;\R^d)})
\le\|\vv\|_{L^2(\varOmega;\R^d)}+\|\nabla\EE(\vv)\|_{L^p(\varOmega;\R^{d\times d\times d})}$.
Then, summing \eqref{test-FF} and \eqref{test-vv} and using 
Gronwall's inequality, we obtain the a-priori estimates
\begin{subequations}\label{est+}\begin{align}\label{est+1}
&\|\vv_{\eps k}\|_{L^\infty(I;L^2(\varOmega;\R^d))\,\cap\,
L^p(I;W_0^{2,p}(\varOmega;\R^d))}^{}\le C\,,
\\[-.3em]&
\|\FF_{\eps k}\|_{L^\infty(I;L^2(\varOmega;\R^{d\times d}))}^{}\le C
\ \ \text{ and }\ \ 
\|\nabla\FF_{\eps k}\|_{L^2(I{\times}\varOmega;;\R^{d\times d\times d}))}^{}\le\frac C{\sqrt\eps}
\,.
\end{align}\end{subequations}

Although these estimates would already allow for passage with $k\to\infty$,
we will still make another a-priori estimate by testing \eqref{Euler2++} in
its Galerkin approximation by $\Delta\FF_{\eps k}$. Using the Green
formula with the boundary condition $(\nn{\cdot}\nabla)\FF_{\eps k}=0$,
we obtain
\begin{align}\nonumber
&\eps\int_\varOmega|\Delta\FF_{\eps k}|^2\,\d x+
\frac{\d}{\d t}\int_\varOmega\frac12|\nabla\FF_{\eps k}|^2\,\d x
=\int_\varOmega\nabla\big((\nabla\vv_{\eps k})\FF_{\eps k}-(\vv_{\eps k}{\cdot}\nabla)\FF_{\eps k}\big)\Vdots\nabla\FF_{\eps k}\,\d x
\\[-.0em]\nonumber&\qquad=
\int_\varOmega(\nabla\FF_{\eps k}{\otimes}\nabla\FF_{\eps k}){:}\EE(\vv_{\eps k})-\frac12|\nabla\FF_{\eps k}|^2{\rm div}\,\vv_{\eps k}
-(\nabla\vv_{\eps k})\nabla\FF_{\eps k}\Vdots\nabla\FF_{\eps k}
-(\nabla^2\vv_{\eps k})\FF_{\eps k}\Vdots\nabla\FF_{\eps k}\,\d x
\\[-.3em]&\qquad\qquad\nonumber
\le\frac52\|\nabla\vv_{\eps k}\|_{L^\infty(\varOmega;\R^{d\times d})}^{}
\|\FF_{\eps k}\|_{L^2(\varOmega;\R^{d\times d})}^2
\\[-.1em]&\hspace{10em}
+\|\nabla^2\vv_{\eps k}\|_{L^p(\varOmega;\R^{d\times d\times d})}^{}
\|\FF_{\eps k}\|_{L^{2^*}(\varOmega;\R^{d\times d})}^{}\big(1+
\|\nabla\FF_{\eps k}\|_{L^2(\varOmega;\R^{d\times d\times d})}^2\big)\,;
\label{test-Delta-FF}\end{align}
here $\nabla\FF{\otimes}\nabla\FF$ denotes the
symmetric matrix $[\nabla\FF{\otimes}\nabla\FF]_{ij}^{}=
\sum_{k,l=1}^d\frac{\partial}{\partial\xx_i}\FF_{kl}^{}
\frac{\partial}{\partial\xx_j}\FF_{kl}^{}$. In \eqref{test-Delta-FF}, we used
both $p\ge d$ so that $p^{-1}+(2^*)^{-1}+2^{-1}\le1$ and also the calculus
\begin{align}\nonumber
&\int_\varOmega\nabla\big((\vv{\cdot}\nabla)\FF
  \big)\Vdots\nabla\FF\,\d x
  =\int_\varOmega(\nabla\FF{\otimes}\nabla\FF){:}\EE(\vv)
  +(\vv{\cdot}\nabla)\nabla\FF\Vdots\nabla\FF\,\d x
\\[-.4em]&\ \ \nonumber
=\int_\varGamma|\nabla\FF|^2\vv{\cdot}\nn\,\d S
+\!\int_\varOmega(\nabla\FF{\otimes}\nabla\FF){:}\EE(\vv)-
({\rm div}\,\vv)|\nabla\FF|^2-\nabla\FF\Vdots
(\vv{\cdot}\nabla)\nabla\FF\,\d x
\\[-.4em]&\hspace{2em}
=\int_\varGamma\frac{|\nabla\FF|^2\!\!}2\ \vv{\cdot}\nn\,\d S
+\int_\varOmega(\nabla\FF{\otimes}\nabla\FF){:}\EE(\vv)-({\rm div}\,\vv)\frac{|\nabla\FF|^2\!\!}2\ \d x\,.
\label{test-Delta}\end{align}
By Gronwall's inequality applied to
\eqref{test-Delta-FF} together with the qualification of the approximate
initial condition $\FF_{0,k}$ bounded in $H^1(\varOmega;\R^d)$,
we obtain still the a-priori estimates
\begin{align}\label{est++}
&\|\FF_{\eps k}\|_{L^\infty(I;H^1(\varOmega;\R^{d\times d}))}^{}\le C
\ \ \ \text{ and }\ \ \ 
\|\Delta\FF_{\eps k}\|_{L^2(I{\times}\varOmega;\R^{d\times d}))}^{}\le\frac C{\sqrt\eps}\,.
\end{align}

\medskip{\it Step 3: limit passage with $k\to\infty$ and with $\eps\to0$}.
Exploiting the estimates \eqref{est+1} and \eqref{est++}, we can select a
subsequence converging in the weak* topologies
$(L^\infty(I;L^2(\varOmega;\R^d))\,\cap\,
L^p(I;W_0^{2,p}(\varOmega;\R^d))\times L^\infty(I;H^1(\varOmega;\R^{d\times d}))$.
The convergence in the
Galerkin approximation is then routine. The highest-order nonlinear term
${\rm div}^2(\nu|\nabla\EE(\vv)|^{p-2}\nabla\EE(\vv))$
is monotone and can be handled by Minty's trick,
while the other terms are linear or, in the case of $(\vv{\cdot}\nabla)\vv$,
$({\rm div}\,\vv)\vv$, $(\vv{\cdot}\nabla)\FF$, $\varphi'(F)F^\top$
and $\varphi(F)$, nonlinear but of a lower order and thus to be
handled by compactness. Here the Aubin-Lions
compact-embedding theorem can be used when employing
estimates of $\frac{\partial}{\partial t}\vv_{\eps k}$ and
$\frac{\partial}{\partial t}\FF_{\eps k}$ by comparison;
cf.\ \cite[Sect.\,8.4]{Roub13NPDE} for adaptations of
that theorem for Galerkin method. In fact, the uniform monotonicity of
$\vv\mapsto{\rm div}^2(\nu|\nabla\EE(\vv)|^{p-2}\nabla\EE(\vv))$,
the weak* convergence of $\vv_{\eps k}$ can be improved to the
strong convergence, so that the Minty trick and the Aubin-Lions
theorem and the estimate of $\frac{\partial}{\partial t}\vv_{\eps k}$
are actually not needed. Also $\frac{\partial}{\partial t}\FF_{\eps k}
=(\nabla\vv_{\eps k})\FF_{\eps k}-(\vv_{\eps k}{\cdot}\nabla)\FF_{\eps k}$ can
be estimated bounded in $L^p(I;L^2(\varOmega;\R^{d\times d}))$ due to \eqref{est+1}
and the first estimate in \eqref{est++}, so that the conventional version of
the Aubin-Lions theorem can be used.

The regularizing term $\eps\Delta\FF_{\eps k}$ vanishes in the limit for
$\eps\to0$ because
$\|\eps\Delta\FF_{\eps k}\|_{L^2(I{\times}\varOmega;\R^{d\times d})}^{}
=\mathscr{O}(\sqrt\eps)\to0$ or, in the weak formulation even 
$\|\eps\nabla\FF_{\eps k}\Vdots\nabla\widetilde\FF\|_{L^1(I{\times}\varOmega)}^{}
=\mathscr{O}(\eps)\to0$. 
Actually, we can pass to the limit simultaneously with $k\to\infty$ and with
$\eps\to0$ since we do not rely on the latter estimate in \eqref{est+1}.
Thus, we obtain a weak solution to the original problem
\eqref{Euler+}.

\medskip{\it Step 4: test by $\SS=\varphi'(\FF)$}.
For the continuous problem \eqref{Euler1+}, we can
perform the physically relevant test by $\vv$ and $\varphi'(\FF)$;
recall that now we have already $\eps=0$. This gives the energy balance
\eqref{energy+++}.
Here it is important that $\pdt{}\vv\in L^p(I;W_0^{2,p}(\varOmega;\R^d)^*)+
L^1(I;L^2(\varOmega;\R^d))$ is in duality with
$\vv\in L^p(I;W_0^{2,p}(\varOmega;\R^d))\cap
L^\infty(I;L^2(\varOmega;\R^d))$, so that the test of \eqref{Euler1+}
by $\vv$ is legitimate. Also, as in Step~3, we have $\pdt{}\FF
=(\nabla\vv)\FF-(\vv{\cdot}\nabla)\FF\in L^p(I;L^2(\varOmega;\R^{d\times d}))$ is
surely in duality
with $\SS=\varphi'(\FF)\in L^\infty(I{\times}\varOmega;\R^{d\times d})$, so that
the physical test of \eqref{Euler2+} by $\SS$ is legitimate. In particular,
all the integrands in \eqref{Sstr-calculus} are integrable (some of them
even bounded). Thus all the calculations leading to the
energy balance  \eqref{energy+++} are not formal.

Obviously, $\nabla\SS=\nabla\varphi'(\FF)
=\varphi''(\FF)\nabla\FF\in L^\infty(I;L^2(\varOmega;\R^{d\times d\times d}))$
provided $\varphi''$ is bounded.
\end{proof}

\begin{remark}[{\sl Stored energy complying with \eqref{ass1}}]\label{example}\upshape
An example of a mechanically relevant frame-indifferent stored energy is
$$
\varphi(F)=\phi(E)\ \ \ \text{ with }\ \
\phi(E)=
\frac{dK|{\rm sph}\,E|^2}{2{+}\eta|E|^{3/2}}
+\frac{G|{\rm dev}\,E|^2}{1{+}\eta|E|^{3/2}}
\ \ \ \text{ where }\ \ E=E(F)=\frac12(F^\top\!F-\bbI)
$$
with the spherical and the deviatoric parts ${\rm sph}\,E=({\rm tr}\,E)\bbI/d$
and ${\rm dev}\,E=E-({\rm tr}\,E)\bbI/d$, and with $K$ and $G$ the bulk and
the shear elastic moduli, respectively. The philosophy of this model is a
quadratic function (here up to higher-order terms in the neighbourhood of 0) of
the Green-Lagrange (sometimes called Green-St.\,Venant) strain tensor $E$.
The regularization parameter  $\eta>0$  is expectedly small just to
ensure that $\varphi(F)=\mathscr{O}(|F|)$ and
$\varphi'(F)=\phi'(F^\top\!F{-}\bbI)E'(F)=\mathscr{O}(1)$ for $|F|\to\infty$,
so that this $\varphi$ complies with \eqref{ass1}. For $\eta=0$
we obtain the
{\it isotropic  St.\,Venant-Kirchhoff material} for which $\varphi$ has,
however, the growth $\mathscr{O}(|F|^4)$ and thus does not comply with
\eqref{ass1}.
\end{remark}

\begin{remark}[{\sl An incompressible limit}]\label{rem-incompressible}\upshape
In literature, the model of the type \eqref{Euler+} is sometimes interpreted
rather as  a  viscoelastic fluid than solid, and then considered as
incompressible.
This is motivated by  a  qualitative difference of bulk and shear
elastic moduli in fluids (the  latter one being  zero) in contrast to
solids where these moduli are mostly of the same order (except rubber-like
materials).
The incompressible model thus modifies \eqref{Euler+} as
\begin{subequations}\label{Euler-incompressible}\begin{align}\nonumber
&     \varrho\DT\vv={\rm div}\big(\bm{T}+\bm{D}\big)+\ff\ \ \text{ and }
\ \ {\rm div}\,\vv=0\,,\ \  \text{where }\ 
\bm{T}=\varphi'(\FF)\FF^\top\!\!+\pi\bbI    \\[-.3em]
    &\hspace*{18.3em}
 \text{ and }\ \ \ \,\bm{D}=\bbD\EE(\vv)-{\rm div}\big(\nu|\nabla \EE(\vv)|^{p-2}\nabla\EE(\vv)\big)\,,
\\[-.3em]\label{Euler2++}
&\DT\FF=(\nabla\vv)\FF \end{align}\end{subequations}
with $\pi$ a pressure.
Sometimes, \eqref{Euler2++} is considered with the nonlinear holonomic
constraint $\det\FF=1$ while  the linear constraint 
${\rm div}\,\vv=0$ is possibly omitted, relying
on that ${\rm div}\,\vv=0$ is equivalent with $\det\FF=1$ if $\det\FF_0=1$, 
cf.\ e.g.\ \cite{CheZha06GESS,HuLin16GSTD,LiLiZh05HVF,Ruzi92MPTM}.
The weak formulation modifies Definition~\ref{def} by imposing
${\rm div}\,\vv=0$ and ${\rm div}\,\widetilde\vv=0$, i.e.\
in particular \eqref{def1} omits the terms
$\frac\varrho2({\rm div}\,\vv)\vv{\cdot}\widetilde\vv$
and $\varphi(\FF)\widetilde\vv$ and \eqref{def2} omits
$({\rm div}\,\vv)\FF{:}\widetilde\FF$. This model naturally arises as
the limit of \eqref{Euler+} by two ways. For the isotropic viscosity tensor
$\bbD\EE(\vv)=K{\rm div}\,\vv+2G{\rm dev}\,\EE(\vv)$ with
${\rm dev}\,\EE=\EE{-}({\rm tr}\,\EE)\bbI/d$ denoting the deviatoric strain
rate, $K$ the bulk modulus, and $G$ the shear modulus, 
the incompressible limit can arise when sending $K\to\infty$. Indeed,
one can estimate
$\|{\rm div}\,\vv\|_{L^2(I\times\varOmega;\R^d)}^{}=\mathscr{O}(1/\sqrt{K})$.
The other way is to assume $\varphi(F)\ge K(1{-}\det F)^2$, which leads
to $\|1{-}\det\FF\|_{L^\infty(I;L^2(\varOmega))}^{}=\mathscr{O}(1/\sqrt{K})$. In both
cases, the limit passage to the weak solutions of the incompressible system
\eqref{Euler-incompressible} is quite straightforward.
\end{remark}

\begin{remark}[{\sl Local non-interpenetration}]\label{rem-blow-up}\upshape
Another physically relevant assumption beside frame indifference is that $\FF$
ranges only ${\rm GL}^+(d)=\{F\in\R^{d\times d};\ \det\,F>0\}$, i.e.\ the
subgroup of the general linear group of matrices with positive determinants.
In particular, one should impose the blow-up assumption $\varphi(F)\to\infty$ 
for $\det F\to0+$, which is however not compatible with \eqref{ass1}.
In some context, one can assume $\varphi:\R^{d\times d}
\to[0,+\infty]$ continuously differentiable on ${\rm GL}^+(d)$ and  
\begin{align}
&\exists\,\epsilon>0, \ \forall F\in\R^{d\times d}:\ \
\varphi(F)\ge\begin{cases}\epsilon
/(\det F)^r\!\!\!&\text{if }\ \det F>0,
\label{ass-blow-up}
\\[-.2em]\quad+\infty&\text{if }\ \det F\le0,\end{cases}\ \ \ \ 
\end{align}
for $r>qd/(q{-}d)$ if one could ensure that $\FF(t)$ ranges a bounded set in
$W^{1,q}(\varOmega;\R^{d\times d})$ for some $q>d$. Then one could use
Healey-Kr\"omer's arguments \cite{HeaKro09IWSS} to ensure $\det\FF$ away from
zero, like  it is  possible in Lagrangian formulation in
\cite{KruRou19MMCM,MieRou16RIEF,MieRou20TKVR}. Here, formally one could
strengthen \eqref{ass3} for $\FF_0\in W^{1,q}(\varOmega;\R^{d\times d})$ and then
test \eqref{Euler2+} by ${\rm div}(|\nabla\FF|^{q-2}\nabla\FF)$,
which gives an estimate
of $\FF$ in $L^\infty(I;W^{1,q}(\varOmega;\R^{d\times d}))$. Assuming
$p\ge2^*q/(2^*{-}q)$ so that $p^{-1}+(2^*)^{-1}+(q')^{-1}\le1$, 
by the H\"older and Young inequalities, this test  modifies
\eqref{test-Delta-FF} as
\begin{align}\nonumber
&\frac{\d}{\d t}\int_\varOmega\frac1q|\nabla\FF|^q\,\d x
=\int_\varOmega\nabla\big((\vv{\cdot}\nabla)\FF
  -(\nabla\vv)\FF\big)\Vdots|\nabla\FF|^{q-2}\nabla\FF\,\d x
\\[-.0em]\nonumber&=
\int_\varOmega|\nabla\FF|^{q-2}(\nabla\FF{\otimes}\nabla\FF){:}\EE(\vv)-\frac1q|\nabla\FF|^q{\rm div}\,\vv-(\nabla\vv)\nabla\FF\Vdots|\nabla\FF|^{q-2}\nabla\FF
-(\nabla^2\vv)\FF\Vdots|\nabla\FF|^{q-2}\nabla\FF\,\d x
\\&\nonumber
\le\frac{\!2q{+}1\!}q\|\nabla\vv\|_{L^\infty(\varOmega;\R^{d\times d})}^{}
\|\nabla\FF\|_{L^q(\varOmega;\R^{d\times d\times d})}^q\!
+\|\nabla^2\vv\|_{L^p(\varOmega;\R^{d\times d\times d})}^{}
\|\FF\|_{L^{2^*}(\varOmega;\R^{d\times d})}^{}\big(1{+}
\|\nabla\FF\|_{L^{q}(\varOmega;\R^{d\times d\times d})}^{q}\big),
\end{align}
where we used a modification of the calculus \eqref{test-Delta}:
\begin{align}\nonumber
&\int_\varOmega\nabla\big((\vv{\cdot}\nabla)\FF
  \big){:}|\nabla\FF|^{q-2}\nabla\FF\,\d x
  =\int_\varOmega|\nabla\FF|^{q-2}(\nabla\FF{\otimes}\nabla\FF){:}\EE(\vv)
  +(\vv{\cdot}\nabla)\nabla\FF\Vdots|\nabla\FF|^{q-2}\nabla\FF\,\d x
\\[-.4em]&\quad\nonumber
=\int_\varGamma|\nabla\FF|^q\vv{\cdot}\nn\,d S
+\int_\varOmega|\nabla\FF|^{q-2}(\nabla\FF{\otimes}\nabla\FF){:}\EE(\vv)-
({\rm div}\,\vv)|\nabla\FF|^q-(q{-}1)|\nabla\FF|^{q-2}\nabla\FF\Vdots
(\vv{\cdot}\nabla)\nabla\FF\,\d x
\\[-.4em]&\hspace{2em}\nonumber
=\int_\varGamma\frac{|\nabla\FF|^q\!\!}q\ \vv{\cdot}\nn\,d S
+\int_\varOmega|\nabla\FF|^{q-2}(\nabla\FF{\otimes}\nabla\FF){:}\EE(\vv)-({\rm div}\,\vv)\frac{|\nabla\FF|^q\!\!}q\ \d x\,.
\end{align}
The boundary integral vanishes if $\vv{\cdot}\nn=0$. Then,
one formally obtains the estimate
$\|\nabla\FF\|_{L^\infty(I;L^q(\varOmega;\R^{d\times d\times d}))}^{}\le C$.
Actually, to  execute this strategy legitimately would be quite technical.
First, a suitable cut-off of the stresses in the momentum equation
is to be done together with strengthening of the $\eps$-regularization
of \eqref{Euler2++} by a considering an $q$-Laplacian
${\rm div}(|\nabla\FF|^{q-2}\nabla\FF)$.
The Galerkin discretization of such a system must be done
separately for the momentum and for the (modified) transport-and-evolution
equation \eqref{Euler2++}, the limit passage of the latter one being performed
first. Also the a-priori estimates are to be extended, in particular
by proving that ${\rm div}(|\nabla\FF|^{q-2}\nabla\FF)$ is in
$L^2(I{\times}\varOmega;\R^{d\times d})$; cf.\ \cite{Roub21QHLS} for such a
strategy in the context of plastic enhancement of the model. Alternatively,
without relying on \eqref{ass-blow-up}, one can
employ directly the transport equation of $1/\det\FF$, namely
$\DT{\overline{1/\det\FF}}=-({\rm div}\,\vv)/\det\FF$.
 To avoid all these very technical arguments,
we confined ourselves to the simpler although less physically relevant
model with \eqref{ass1} instead of \eqref{ass-blow-up}.\end{remark}

\begin{remark}[{\sl Reconstruction of an underlying
deformation}]\label{rem-deformation}\upshape
Implicitly, we have in mind the situation when 
$\FF_0=\nabla_\XX^{}\yy_0$ for some
initial deformation $\yy_0\in H^2(\varOmega;\R^d)$. Although we
have not explicitly needed this additional qualification of $\FF_0$,
in the context of the original motivation of the model, a natural question
is whether one can reconstruct the deformation $\yy(t)$
 such that $\vv=\DT\yy$ and $\FF=\nabla_\XX^{}\yy$, and also $\yy(0)=\yy_0$.
This seems a nontrivial question, however. We can always construct the return
mapping $\bm\xi$ mentioned above by solving the simple transport equation
$\DT{\bm\xi}=\bm0$ with the initial condition $\bm\xi(0)$=identity. Then
$\FF=(\nabla_\xx\bm\xi)^{-1}$ and, if
$\bm\xi(t):\varOmega\to\varOmega$ is surjective, $\yy(t)=\bm\xi^{-1}(t)$.
This global surjectivity seems
not automatic, however. An example for such global surjectivity, indicating
the complexity of this problem,
is for a completely fixed boundary deformation (not considered in
this paper), i.e., in addition to $\vv{\cdot}\nn=0$, also the tangential
velocity on $\varGamma$ would be prescribed zero so that $\vv=\bm0$ on
$\varGamma$. Then also $\bm\xi|_{\varGamma}^{}(t)=\,$identity on $\varGamma$.  As
the velocity field $\vv$ is enough regular, the regularity of the initial
condition $\bm\xi_0\,$=\,identity and local surjectivity in the sense
of invertibility of $\nabla\bm\xi_0$, i.e.\ here
$\det(\nabla\bm\xi_0)=\det(\bbI)=1>0$, is copied for all $t>0$. Note
that we have the evolution-and-transport equation
$\DT{\overline{\det\nabla\bm\xi}}=-({\rm div}\,\vv)\det\nabla\bm\xi$, as
actually mentioned in Remark~\ref{rem-blow-up} since
$\det\nabla\bm\xi=1/\det\FF$. Then the
classical result of J.M.\,Ball \cite{Ball81GISF} shows global injectivity
of $\bm\xi(t)$, i.e.\ $\bm\yy(t)=\bm\xi^{-1}(t)$ exists.
\end{remark}

\section*{Acknowledgments}
The author is thankful to Ulisse Stefanelli and Giuseppe Tomassetti for
valuable discussions about the Eulerian continuum mechanics.

\end{sloppypar}
\end{document}